\documentclass[a4paper,11pt]{article}

\usepackage[T1]{fontenc}
\usepackage[utf8]{inputenc}
\usepackage[english]{babel}
\usepackage{geometry}
\usepackage{amsthm}
\usepackage{amsmath}
\usepackage{graphicx}
\usepackage{amssymb}
\usepackage[noadjust]{cite}
\usepackage{indentfirst}
\usepackage{hyperref}

\usepackage{enumitem}
\setlist[enumerate]{font=\upshape,label=(\arabic*)}
\setlist{itemsep=1pt,topsep=4pt}

\numberwithin{equation}{section}
\allowdisplaybreaks

\usepackage{tikz}
\usetikzlibrary{arrows.meta}

\newtheorem{statement}{Statement}[section]
\newtheorem{theorem}[statement]{Theorem}
\newtheorem{lemma}[statement]{Lemma}
\newtheorem{corollary}[statement]{Corollary}
\newtheorem{proposition}[statement]{Proposition}
\theoremstyle{definition}
\newtheorem{definition}[statement]{Definition}
\newtheorem{example}[statement]{Example}
\newtheorem{notation}[statement]{Notation}
\newtheorem{remark}[statement]{Remark}

\newcommand{\pp}{\mathbb{P}}
\newcommand{\pn}{\mathbb{P}_n}
\newcommand{\ppnm}{\mathbb{P}_{n\times m}}

\newcommand{\crk}{\mathrm{crk}}
\newcommand{\ci}{\mathrm{CI}}

\newcommand{\divides}{%
  \mathrel{%
    \nonscript\mspace{-\thickmuskip}%
    \nonscript\mspace{0.5\thinmuskip}%
    |
    \nonscript\mspace{0.5\thinmuskip}%
    \nonscript\mspace{-\thickmuskip}%
  }%
}

\providecommand{\keywords}[1]{\medskip\noindent\textbf{Keywords:} #1}
\providecommand{\msc}[1]{\noindent\textbf{MSC 2020:} #1}

\newcommand\freefootnote[1]{%
    \bgroup
    \renewcommand\thefootnote{\fnsymbol{footnote}}%
    \renewcommand\thempfootnote{\fnsymbol{mpfootnote}}%
    \footnotetext[0]{#1}%
    \egroup
}

\date{}

\title{Graphs with the largest chainable index\\ and their properties}
\author{Alexander Guterman$^{a}$,
Svetlana Zhilina$^{b,c}$,
Evgenii Shafeev$^{b,c}$
}
\date{\small \em
$^a$Bar-Ilan University, Ramat Gan, 5290002, Israel\\
$^b$Lomonosov Moscow State University, Moscow, 119991, Russia\\
$^c$Moscow Center of Fundamental and Applied Mathematics, Moscow, 119991, Russia
}
\begin{document}
\maketitle

\begin{abstract}
This paper provides a complete characterization of the set of graphs with the largest chainable index. In particular, it is shown that the number of arcs and vertex degrees of such graphs satisfy strong restrictions. Moreover, graphs with the largest chainable index form five parametric families, and exactly one of them consists of graphs that are not strongly connected. As a consequence, the number of graphs on $n$ vertices with the largest chainable index is computed, and the values of the primitive exponent for these graphs are calculated.
\end{abstract}


\msc{05C50, 05C35, 05C38.}

\freefootnote{The work of the second and the third authors was supported by the Ministry of Science and Higher Education of the Russian Federation as part of the program of the Moscow Center of Fundamental and Applied Mathematics under agreement No. 075-15-2025-345.}

\freefootnote{Email addresses: \texttt{alexander.guterman@biu.ac.il} (Alexander Guterman), \texttt{s.a.zhilina@gmail.com} (Svetlana Zhilina), \texttt{shafeev.ev@yandex.ru} (Evgenii Shafeev)}

\section{Introduction}

\begin{notation}
We denote by $\pp$ the set of all real nonnegative matrices without zero rows and columns. By $\ppnm$ we denote the set of $n\times m$ matrices from $\pp$, and by $\pn$ the set $\pp_{n\times n}$.
\end{notation}

\begin{definition}
A nonnegative matrix $A = (a_{ij}) \in \pp$ is called \textit{chainable} if for every pair of its positive entries $a_{ij}$, $a_{pq}$ there exists a chain with these entries as endpoints, that is, a sequence of positive entries $a_{i_1j_1},a_{i_2j_2},\ldots,a_{i_mj_m}$ such that:
\begin{enumerate}
    \item $(i_1, j_1) = (i, j)$;
    \item $(i_m, j_m) = (p, q)$;
    \item for each $l = 1,2,\ldots,m-1$ we have $i_l = i_{l+1}$ or $j_l = j_{l+1}$.
\end{enumerate}
\end{definition}

If we imagine the entries of a matrix as the squares of a chessboard and a rook is allowed to move only from one positive entry to another, then the definition can be rephrased as follows. A matrix is chainable if the rook can get from any positive square to any other one.

Chainable matrices were first considered in~\cite{Hartfiel} and~\cite{Sinkhorn}. In the present paper, when speaking about square matrices of order $n$, we always assume that $n \geq 2$.

\begin{definition}
A matrix $A \in \pn$ is called \textit{potentially chainable} if for some $k \in \mathbb{N}$ the matrix $A^k$ is chainable. The smallest such~$k$ is called the \textit{chainable index}. The chainable index of a matrix $A$ is denoted by $\ci(A)$. If $A$ is not potentially chainable, then we set $\ci(A) = 0$.
\end{definition}

The notions of a potentially chainable matrix and of the chainable index were introduced and studied by Alpin and Bashkin in~\cite{AlpinBashkin}. This research was continued, in particular, in~\cite{AlpinGutermanShaf,Frob,preservers}.

\begin{theorem}[{\cite[Theorem 3.8]{AlpinGutermanShaf}}]
Let $A \in \pn$. Then $\ci(A) \leq n-1$.
\end{theorem}

\begin{theorem}[{\cite[Corollary 4.8]{AlpinGutermanShaf}}]
All values of the chainable index from $0$ to $n-1$ are attained by matrices in~$\pn$.
\end{theorem}

The chainable index is related to another important invariant of nonnegative matrices, namely, the primitive exponent.

\begin{definition}
A nonnegative matrix $A \in \pn$ is called {\em primitive} if for some $k \in \mathbb{N}$ the matrix $A^k$ is (entrywise) positive. The smallest such $k$ is called the {\em exponent} (the {\em primitive exponent}) of $A$ and is denoted by $\exp(A)$. If $A$ is not primitive, then we set $\exp(A) = 0$.
\end{definition}

\begin{definition}
Two nonnegative matrices $A = (a_{ij})$ and $B = (b_{ij})$ of the same size are called {\em combinatorially equivalent} (denoted by $A \sim B$) if $a_{ij} > 0$ if and only if $b_{ij}>0$.
\end{definition}

It is clear from the definitions that combinatorially equivalent matrices have equal exponents. The same is true for the chainable index.

\begin{definition}
The Wielandt matrix $W_n$, $n \geq 2$, is defined as follows:
\[W_n = \begin{pmatrix}
    0 & 1 & 0 & \dots & 0\\
    0 & 0 & 1 & \dots & 0\\
    \dots & \dots & \dots & \dots & \dots \\
    0 & 0 & 0 & \dots & 1\\
    1 & 1 & 0 & \dots & 0
\end{pmatrix}.\]
\end{definition}

\begin{theorem}[{\cite{Wielandt}}]\label{THexponent}
Let $A$ be a nonnegative $n\times n$ matrix. Then
\begin{enumerate}
    \item $\exp(A) \leq (n-1)^2 + 1$;
    \item $\exp(A) = (n-1)^2 + 1$ if and only if there exists a permutation matrix $P$ such that $PAP^T \sim W_n$.
\end{enumerate}
\end{theorem}

\begin{theorem}[{\cite[p.~23]{AlpinBashkin}}]
$\ci(W_n) = n-1$.
\end{theorem}
Thus $W_n$ attains the largest possible value of the chainable index; in what follows we say simply that it has the largest chainable index. Theorem~\ref{THexponent} gives a complete description of the set of matrices with the largest exponent. However, the following example shows that a matrix with the largest chainable index need not have the largest exponent.

\begin{example}\label{kexample}
Consider the matrix
\[M = \begin{pmatrix}
1 & 1 & 0 & 0 \\
0 & 0 & 1 & 0 \\
0 & 0 & 0 & 1 \\
0 & 0 & 0 & 1 \\
\end{pmatrix},
\quad M^2 = \begin{pmatrix}
1 & 1 & 1 & 0 \\
0 & 0 & 0 & 1 \\
0 & 0 & 0 & 1 \\
0 & 0 & 0 & 1 \\
\end{pmatrix},
\quad M^k \sim \begin{pmatrix}
1 & 1 & 1 & 1 \\
0 & 0 & 0 & 1 \\
0 & 0 & 0 & 1 \\
0 & 0 & 0 & 1 \\
\end{pmatrix}\]
for all $k \geq 3$. The matrix $M^3$ is chainable, while $M^2$ is not, since the first row of $M^2$ does not intersect any other row. Therefore, $\ci(M) = 3$. At the same time, the matrix $M$ is not primitive, whereas $W_4$ is primitive. Consequently, there is no permutation matrix $P$ such that $PMP^T \sim W_4$. Thus, there exist matrices with the largest chainable index that are not combinatorially equivalent to $W_4$ under any simultaneous permutation of rows and columns.
\end{example}

\begin{definition}
A square nonnegative matrix $A$ is called {\em irreducible} if for any $i,j$ there exists $k \in \mathbb{N}$ such that $(A^k)_{ij} > 0$.
\end{definition}

It immediately follows from the definition that any primitive matrix is irreducible.

\begin{proposition}[{\cite[Ch.~13, Lemma~1]{Gantmacher}}]
    Let $A$ be a square nonnegative matrix of order~$n$. The following conditions are equivalent:
    \begin{enumerate}
        \item $A$ is irreducible;
        \item $(I+A)^{n-1} > 0$;
        \item there is no permutation matrix $P$ such that
        \[PAP^T = \begin{pmatrix}
            A_1 & A_2 \\
            0 & A_3\\
        \end{pmatrix},\]
        where $A_1, A_3$ are nonempty square blocks.
    \end{enumerate}
\end{proposition}

It can be easily seen that every primitive matrix is potentially chainable; however, by Example~\ref{kexample}, the converse is false. Nevertheless, as the following theorem shows, in the case of irreducible matrices these two notions are equivalent. This means that the class of potentially chainable matrices is a natural extension of the class of primitive matrices.

\begin{theorem}[{\cite[Proposition 6.1]{AlpinBashkin}}]\label{theoremprimtochain}
Let $A$ be an irreducible matrix. Then $A$ is potentially chainable if and only if it is primitive.
\end{theorem}

When analyzing combinatorial properties of nonnegative matrices, one often uses a graph-theoretic interpretation. In this paper we consider directed graphs without multiple arcs, but possibly with loops. Moreover, we assume that a graph has a finite number of vertices, which is at least two.

To each nonnegative $n \times n$ matrix $A$ we associate a directed graph $D = D(A)$ with the vertex set $V(D) = \{1, 2, \ldots, n\}$ and the arc set $E(D) = \{(i,j) \mid a_{ij} > 0\}$. Thus, the matrix $A$ is combinatorially equivalent to the adjacency matrix of the graph~$D$.

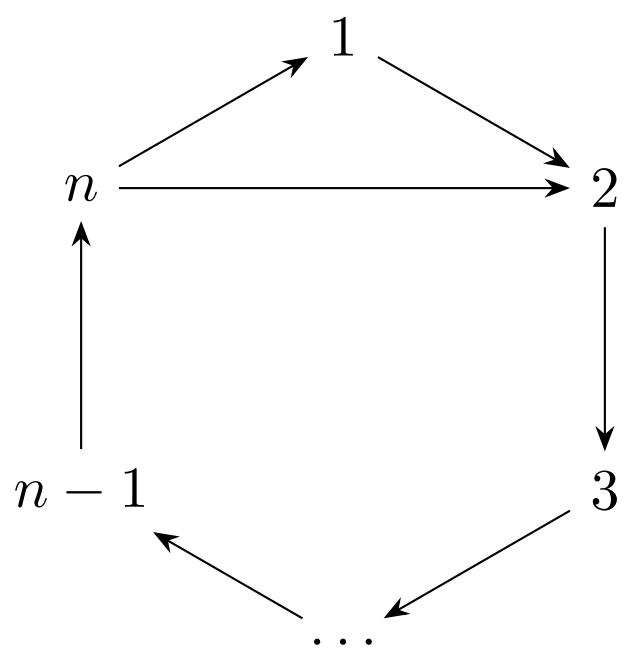
\begin{figure}[ht]       
    \centering
    \begin{tikzpicture}[
    > = Stealth,          
    vertex/.style = {
        font=\sffamily
    }]
    \def\R{2.0}
    \node (1) at (90:\R)  {$1$};
    \node (2) at (30:\R)  {$2$};
    \node (3) at (-30:\R) {$3$};
    \node (4) at (-90:\R) {$\ldots$};
    \node (5) at (210:\R) {$n-1$};
    \node (6) at (150:\R) {$n$};

    \draw[->] (1) -- (2);
    \draw[->] (2) -- (3);
    \draw[->] (3) -- (4);
    \draw[->] (4) -- (5);
    \draw[->] (5) -- (6);
    \draw[->] (6) -- (1);
    \draw[->] (6) -- (2);
\end{tikzpicture}
 \caption{The Wielandt graph $W_n$, $n \geq 3$.}
 \label{fig:Wielandt}
\end{figure}

As an example, consider the graph of the Wielandt matrix for $n \geq 3$, shown in Fig.~\ref{fig:Wielandt}. We use the notation $W_n$ both for the Wielandt matrix and for the corresponding graph.

Recall that a walk may repeat vertices and arcs, whereas a path is a walk without repeated vertices. We remark that a shortest walk between two vertices is always a path.

\begin{notation}
    The notation $u\xrightarrow{t} v$ means that there exists a directed walk from $u$ to~$v$ of length exactly $t$. We write $u\rightarrow v$ when $u \xrightarrow{1} v$. For convenience, we assume that $u \xrightarrow{0} u$ for any vertex~$u$.
\end{notation}

\begin{definition}
A directed graph is called {\em strongly connected} if for any pair of vertices $u, v$ we have $u \xrightarrow{t} v$ for some $t \geq 0$. A directed graph is called {\em (weakly) connected} if, ignoring the directions of arcs, one can get from any vertex $u$ to any vertex~$v$. Thus, any strongly connected graph is connected.
\end{definition}

For example, the graph $W_n$ is strongly connected, but becomes merely connected after the removal of any arc except~$(n,2)$.

The following facts on the correspondence between matrices and graphs are well known and follow from the definition of the graph $D(A)$.
\begin{proposition}\label{propbasic}
Let $A$ be a nonnegative matrix.
\begin{enumerate}
    \item The graph $D^T = D(A^T)$ is obtained from $D(A)$ by reversing the directions of all arcs.
    \item Let $k \in \mathbb{N}_0$. Then $u \rightarrow v$ in the graph $D^k = D(A^k)$ if and only if there is a walk $u \xrightarrow{k} v$ in $D(A)$.
    \item The graph $D(A)$ is strongly connected if and only if the matrix $A$ is irreducible.
\end{enumerate}
\end{proposition}

We extend any property or numerical invariant defined for a matrix $A \in\pn$ to the graph $D(A)$. In particular, the chainable index of a graph is the chainable index of its adjacency matrix. In the same way we define the exponent of a graph, as well as the notions of a primitive and a potentially chainable graph.
For example, a graph $D$ is primitive if for some $t \in \mathbb{N}$ we have $i\xrightarrow{t} j$ for all $i,j \in V(D)$. The smallest such $t$ is the exponent of the graph~$D$.

From now on we consider only matrices that belong to the set $\pp$. Similarly, among directed graphs on $n$ vertices we consider only those whose adjacency matrix belongs to~$\pp_n$. Thus, the chainable index of such graphs is well defined.

\medskip

This paper classifies all graphs which attain the largest chainable index and investigates their properties. We show that, in such a graph, at most four vertices have more than two incident arcs. Any other vertex is incident to exactly two distinct arcs: one incoming and one outgoing. These strong constraints on vertex degrees allow us to describe the structure of all such graphs explicitly.

Graphs on $n$ vertices that attain the largest chainable index can have either $n+1$ or $n+2$ arcs. These two types of graphs are considered separately. The first type splits into three families of graphs. Only one of these families consists of graphs that are not strongly connected; its members are known as dumbbell graphs. Despite this difference, the criterion for the chainable index to be the largest is the same for all families: the GCD of the set of lengths of all simple directed cycles must be equal to one. The other two families are the well-known bicyclic graphs $\Theta$ and~$\infty$. The Wielandt graph~$W_n$, whose exponent is the largest, is a special case of the family~$\Theta$. The second type contains two families of graphs, both consisting of strongly connected graphs. The first one consists of a single graph $W_n'$, which differs from the graph $W_n$ by the addition of a loop $1\rightarrow 1$; the second family $\overline{\infty}$ is structurally similar to the family~$\infty$. The structure of all graphs with the largest chainable index is simple enough for their exponents to be expressed explicitly in terms of the lengths of their independent cycles.

The paper is organized as follows. Section~$2$ introduces the notion of the chainable rank and describes the necessary properties of graphs with the largest chainable index. In particular, Lemma~\ref{219} shows that these graphs satisfy certain restrictions on their vertex degrees and fall into two types, depending on the total number of arcs. In Section~$3$ we obtain a characterization of graphs of the first type, and in Section~$4$ of graphs of the second type, see Theorems~\ref{theorem:typeoneclassification} and~\ref{theorem:typetwoclassification}, respectively. As a result, the set of graphs with the largest chainable index is completely described, see Theorem~\ref{theorem:fullclassification}. Section~$5$ is devoted to computing the exponents of these graphs, and Section~$6$ determines their number for a fixed number of vertices.

\section{Properties of graphs with the largest chainable index}

The following notions are needed for a more detailed analysis of potentially chainable matrices and graphs.

\begin{definition}\label{definition:rowintersect}
Let $A = (a_{il}) \in \pp_{n\times m}$. We say that the rows $i$ and~$j$ of the matrix $A$ {\em intersect} (denoted by $i \wedge j$) if there exists $k \in \{ 1, \ldots, m \}$ such that $a_{ik} > 0$ and $a_{jk} > 0$.
\end{definition}

\begin{definition}\label{chainlink}
Let $A \in \ppnm$. We say that the rows $i$ and~$j$ are {\em linked by a chain} if there exists a sequence of rows $k_1, \ldots, k_t$ such that
\[
(i = k_1) \wedge k_2 \wedge \dots \wedge k_{t-1} \wedge (k_t = j),
\]
that is, $k_l$ intersects $k_{l+1}$ for $l = 1, \ldots, t - 1$.
\end{definition}

\begin{proposition}[{\cite[Lemma 2.1]{AlpinBashkin}}]
For a matrix $A \in \ppnm$, the relation of being linked by a chain is an equivalence relation on the set of its rows.
\end{proposition}

\begin{definition}[{\cite[Definition 3.1]{AlpinBashkin}}]\label{Def3.4}
For a matrix $A \in \ppnm$, the number of equivalence classes with respect to the relation of being linked by a chain is called the {\em chainable rank} and is denoted by $\crk(A)$.
\end{definition}
The chainable rank of a matrix $A \in \ppnm$ takes values from $1$ to $\min\{n,m\}$.
\begin{remark}
A matrix $A \in \ppnm$ is chainable if and only if $\crk(A) = 1$.
\end{remark}

The chainable rank has many properties of the usual matrix rank, see~\cite{AlpinBashkin,AlpinGutermanShaf,Frob}. We will use the following fact:

\begin{theorem}[{\cite[Theorem 3.1]{AlpinBashkin}}]\label{thcrkab}
If for matrices $A,B \in \pp$ the product $AB$ exists, then
\[\crk(AB) \leq \min \{\crk(A), \crk(B)\}.\]
\end{theorem}

The chainable rank of a graph is the chainable rank of its adjacency matrix. The following definition is a graph analogue of Definition~\ref{definition:rowintersect}.

\begin{definition}
Vertices $u,v$ of a graph $D$ {\em intersect} at a vertex $w$ after raising to a power $k \in \mathbb{N}$ if $u \xrightarrow{k} w$ and $v \xrightarrow{k} w$. In this case we write $u \wedge^k v$, and we abbreviate $u \wedge^1 v$ as $u \wedge v$.
\end{definition}

\begin{lemma}\label{lemma:CIgraph}
A graph $D$ is potentially chainable if and only if there exists $k\in \mathbb{N}$ such that for any pair of vertices $u,v \in V(D)$ there exists a chain of vertices
\begin{equation} \label{equation:chain}
(u = x_1) \wedge^k x_2 \wedge^k \ldots \wedge^k x_{t-1} \wedge^k (x_t = v),
\end{equation}
that is, $x_l \wedge^k x_{l+1}$ for $l = 1,\ldots, t-1$. The smallest such $k$ is the chainable index of $D$.
\end{lemma}
\begin{proof}
Let $A \in \pn$ be the adjacency matrix of $D$, so $D = D(A)$. Consider any two vertices $u, v \in V(D)$. By Proposition~\ref{propbasic}, we have $u \wedge^k v$ if and only if there exists $w$ such that $(A^k)_{u,w} > 0$ and $(A^k)_{v,w} > 0$. Therefore, condition~\eqref{equation:chain} is equivalent to the matrix $A^k$ being chainable, and the smallest such $k$ coincides with~$\ci(A)$.
\end{proof}
In particular, a graph $D$ is chainable if in the chain above we have $x_l \wedge x_{l+1}$ for all $l = 1,\ldots, t-1$.

\begin{proposition}\label{proposition:pcisweak}
A potentially chainable graph is connected.
\end{proposition}
\begin{proof}
For any pair of vertices $u,v$ in a potentially chainable graph there exists a chain~\eqref{equation:chain}, that is, there are vertices $w_1, \dots, w_{t-1}$ such that
\[(u = x_1) \xrightarrow{k} w_1, \quad x_2 \xrightarrow{k} w_1, \quad \ldots, \quad x_{t-1} \xrightarrow{k} w_{t-1}, \quad (v=x_t) \xrightarrow{k} w_{t-1}.\]
The concatenation of these walks, with all walks of even index reversed, forms an undirected walk connecting $u$ and~$v$.
\end{proof}

\begin{lemma}\label{lemmaonlyPC}
Let $D$ be a potentially chainable graph on $n$ vertices. Then the equality $\crk(D) = n-1$ is equivalent to $\ci(D) = n-1$.
\end{lemma}
\begin{proof}
According to \cite[Corollary 2.13]{preservers}, the chainable rank of a graph with the largest chainable index equals $n-1$. It remains to prove that any potentially chainable graph with chainable rank $n-1$ has the largest chainable index. Recall that for any $k \in \mathbb{N}_0$ we have $D^k = D(A^k)$, where $A$ is the adjacency matrix of~$D$. Then $\crk(D^k) = \crk(A^k)$. We use the inequality
\[\crk\big(D^k\big) - \crk\big(D^{k+1}\big) \leq \crk\big(D^{k-1}\big) - \crk\big(D^k\big), \quad k \in \mathbb{N},\]
obtained in~\cite[Corollary 4.14]{Frob}.
By Theorem~\ref{thcrkab}, the chainable rank of $D^t$ does not increase as $t$ grows. If $\crk(D) = n-1$, then, since $\crk(D^0) = n$, we obtain $\crk(D^{k-1}) - \crk(D^{k}) \leq 1$ for any~$k$. If for some $k < n$ we have $\crk(D^{k-1}) - \crk(D^{k}) = 0$, then for all $m \geq k$ we have $\crk(D^{m-1}) - \crk(D^{m}) = 0$, whence
\begin{align*}
\crk(D^m) = \crk(D^{k-1}) &= \crk(D^0) - (\crk(D^0) - \crk(D^1)) - \ldots\\
&\hphantom{= \crk(D^0) \;} - (\crk(D^{k-2}) - \crk(D^{k-1}))\\
&\geq n - (k - 1) \geq 2.
\end{align*}
That is, the graph $D^m$ is not chainable for all $m \geq k$. This contradicts the fact that $D$ is a potentially chainable graph. Hence, for each $k = 1,\ldots, n-1$ we have $\crk(D^{k-1}) - \crk(D^{k}) = 1$ and, consequently, $\crk(D^k) = n-k$. In particular, $\crk(D^{n-2}) = 2$ and $\crk(D^{n-1}) = 1$, that is, $\ci(D) = n-1$.
\end{proof}

\begin{corollary}[{\cite[Lemma 2.14 and Corollary 2.15]{preservers}}]\label{lemmaonewedge}
Let $D$ be a graph with the largest chainable index, that is, $|V(D)| = n$ and $\ci(D) = n-1$. Then
\begin{enumerate}
    \item there is exactly one pair of distinct vertices $u, v \in V (D)$ such that $u\wedge v$;
    \item there is exactly one pair of distinct vertices $u', v' \in V (D)$ for which there exists a vertex $w$ such that $w \rightarrow u'$, $w \rightarrow v'$.
\end{enumerate}
\end{corollary}

\begin{definition}
The in-degree $d^-(v)$ of a vertex $v \in V(D)$ is the number of arcs entering this vertex. The out-degree $d^+(v)$ is the number of arcs leaving this vertex.
\end{definition}
For example, in the graph $W_n$ in Fig.~\ref{fig:Wielandt} we have $d^+(n) = 2$, $d^+(i) = 1$ for $i \neq n$, and $d^-(2) = 2$, $d^-(j) = 1$ for $j \neq 2$.

\begin{remark}
A matrix $A$ belongs to the set $\pn$ if and only if we have $d^+(v) \geq 1$ and $d^-(v) \geq 1$ for every vertex $v$ in its graph $D(A)$.
\end{remark}

\begin{proposition}\label{propintersect}
If $u \wedge^k v$ in the graph $D(A)$, where $A \in \pn$, then $u \wedge^{k+1} v$.
\end{proposition}
\begin{proof}
Let $u \xrightarrow{k} w$ and $v \xrightarrow{k} w$. Since $A \in \pn$, in the graph $D(A)$ we have $d^+(w) \geq 1$, that is, there is a vertex $w'$ such that $w\rightarrow w'$. Hence $u \xrightarrow{k+1} w'$ and $v\xrightarrow{k+1} w'$.
\end{proof}

\begin{definition} \label{definition:types}
Let $D$ be a graph with the largest chainable index.
\begin{enumerate}
    \item We say that $D$ is of the first type if it has exactly one vertex $z$ with $d^-(z) = 2$, while all other vertices have $d^- = 1$, and exactly one vertex $w$ with $d^+(w) = 2$, while all other vertices have $d^+ = 1$. Here $z$ and $w$ may coincide; in this case the common vertex has $d^+ = d^- = 2$.
    \item We say that $D$ is of the second type if it has exactly two vertices $z_1, z_2$ with $d^-(z_1) = d^-(z_2) = 2$, while all other vertices have $d^- = 1$, and exactly two vertices $w_1, w_2$ with $d^+(w_1) = d^+(w_2) = 2$, while all other vertices have $d^+ = 1$; here $z_i$ may coincide with~$w_j$.
\end{enumerate}
\end{definition}

\begin{lemma}[{\cite[Lemma 2.19]{preservers}}] \label{219}
Let $D$ be a graph with the largest chainable index. Then $D$ is either of the first or of the second type.
\end{lemma}

The strong restrictions on vertex degrees allow us to describe explicitly the structure of graphs with the largest chainable index of both types.

\section{Graphs with the largest chainable index of the first type}

\begin{figure}[htp]
    \centering
    \begin{tikzpicture}[
    > = Stealth,          
    vertex/.style = {
        font=\sffamily
    }]
        \node (1) at (-5,0) {$1$};
        \node (2) at (-6.5,1.5) {$2$};
        \node (3) at (-5,3) {$\dots$};
        \node (4) at (-3.5,1.5) {$a$};
        \node (5) at (-3,0) {$a+1$};
        \node (6) at (-1,0) {$\dots$};
        \node (7) at (1,0)  {$a+b-1$};
        \node (8) at (3,0)  {$a+b$};

        \node (9)  at (1.5,1.5)  {$a+b+1$};
        \node (10) at (3,3)  {$\dots$};
        \node (11) at (4.5,1.5)  {$\qquad a +b + c -1$};

        \draw[->] (1) -- (2);
        \draw[->] (2) -- (3);
        \draw[->] (3) -- (4);
        \draw[->] (4) -- (1);

        \draw[->] (1) -- (5);
        \draw[->] (5) -- (6);
        \draw[->] (6) -- (7);
        \draw[->] (7) -- (8);

        \draw[->] (8) -- (9);
        \draw[->] (9) -- (10);
        \draw[->] (10) -- (11);
        \draw[->] (11) -- (8);

    \end{tikzpicture}
    \caption{The graph $K(a,b,c)$.}
    \label{fig:two_rhombus}
\end{figure}
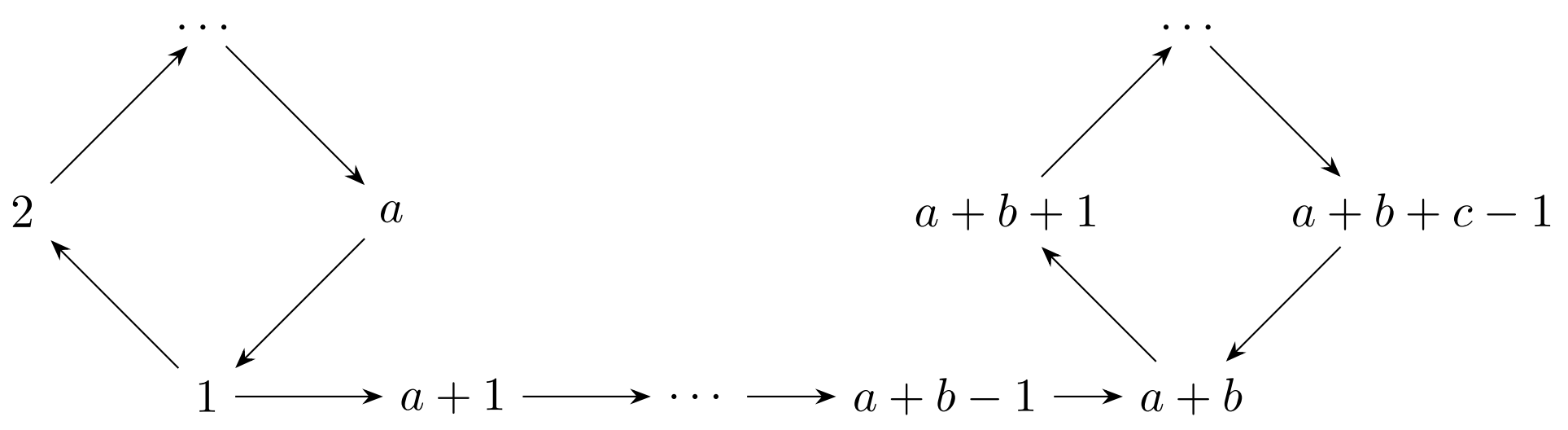

\begin{definition}
Let $a, b, c \in \mathbb{N}$. The graph $K(a,b, c)$, shown in Fig.~\ref{fig:two_rhombus}, is defined as follows:
\begin{gather*}
    1\rightarrow 2\rightarrow \ldots \rightarrow a \rightarrow 1,\\
    1 \rightarrow a+1 \rightarrow \ldots \rightarrow a + b,\\
    a + b \rightarrow a + b +1 \rightarrow \ldots \rightarrow a + b + c-1 \rightarrow a + b.
\end{gather*}
Thus, the graph $K(a,b,c)$ consists of two disjoint simple cycles $\{1,\ldots, a\}$ and $\{a+b, \ldots,\allowbreak a + b + c-1\}$ of lengths $a$ and $c$, respectively, and a simple path of length~$b$ connecting these cycles. If $a = 1$ or $c = 1$, then the corresponding cycle is a loop.
\end{definition}

The matrix $M$ from Example~\ref{kexample} is the adjacency matrix of the graph $K(1,3,1)$.

\begin{definition}
Let $a, b \in \mathbb{N}$, $a+b>2$. The graph $\infty(a,b)$, shown in Fig.~\ref{fig:two_cycles}, is defined as follows:
\begin{gather*}
    1\rightarrow 2 \rightarrow \ldots \rightarrow a \rightarrow 1,\\
    1 \rightarrow a + 1 \rightarrow \ldots \rightarrow a + b-1 \rightarrow 1.
\end{gather*}
Thus, the graph $\infty(a,b)$ consists of two simple cycles of lengths $a$ and~$b$ that intersect in exactly one vertex. If $a = 1$ or $b=1$, then one of these cycles is a loop.
\end{definition}

\begin{definition}
Let $a, b, c \in \mathbb{N}$, $a+b>2$. The graph $\Theta(a,b,c)$, shown in Fig.~\ref{fig:theta}, is defined as follows:
\begin{gather*}
    1 \rightarrow 2 \rightarrow \ldots \rightarrow a+1,\\
    1 \rightarrow a + 2 \rightarrow \ldots \rightarrow a + b \rightarrow a +1,\\
    a +1 \rightarrow a + b+1 \rightarrow\ldots \rightarrow a+b+c -1 \rightarrow 1.
\end{gather*}
If $b = 1$, then the second path has the form $1 \rightarrow a+1$. If $c = 1$, then the third path has the form $a+1 \rightarrow 1$. The graph $\Theta(a,b,c)$ has two cycles of lengths $a+c$ and $b+c$, which intersect along a simple path of length~$c$.
\end{definition}

Its particular case is the Wielandt graph $W_n$, which is isomorphic to~$\Theta(1,2,n-2)$ for $n \geq 3$.

\begin{figure}[t]
\centering
\begin{minipage}[b]{0.37\textwidth}
\centering
\begin{tikzpicture}[
    > = Stealth,          
    vertex/.style = {
        font=\sffamily
    }
]
        \node (1) at (0,0) {$1$};
        \node (2) at (-1.5,1.5) {$2$};   
        \node (3) at (0,3) {$\ldots$};   
        \node (4) at (1.5,1.5) {$a$};  
        \node (5) at (1.5,-1.5) {$a+1$};    
        \node (6) at (0,-3) {$\ldots$};    
        \node (7) at (-1.5,-1.5) {$a + b - 1$};   

        \draw[->] (1) -- (2);
        \draw[->] (2) -- (3);
        \draw[->] (3) -- (4);
        \draw[->] (4) -- (1);

        \draw[->] (1) -- (5);
        \draw[->] (5) -- (6);
        \draw[->] (6) -- (7);
        \draw[->] (7) -- (1);

    \end{tikzpicture}
    \caption{The graph $\infty(a,b)$.}
    \label{fig:two_cycles}
\end{minipage}%
\begin{minipage}[b]{0.63\textwidth}
\centering
\begin{tikzpicture}[
    > = Stealth,          
    vertex/.style = {
        font=\sffamily
    }
]
    \node (1) at (-4,0) {$1$};
    \node (5) at (4,0) {$a+1$};
    \node (2) at (-2.5,1.8) {$2$};
    \node (3) at (0,2.5) {$\ldots$};
    \node (4) at (2.5,1.8) {$a$};
    \node (6) at (-2,0) {$a+2$};
    \node (7) at (0,0) {$\dots$};
    \node (8) at (2,0) {$a+b$};
    \node (9) at (2.5,-1.8) {$a+b+1$};
    \node (10) at (0,-2.5) {$\dots$};
    \node (11) at (-2.5,-1.8) {$a+b+c-1$};

    \draw[->] (1) -- (2);
    \draw[->] (2) -- (3);
    \draw[->] (3) -- (4);
    \draw[->] (4) -- (5);

    \draw[->] (1) -- (6);
    \draw[->] (6) -- (7);
    \draw[->] (7) -- (8);
    \draw[->] (8) -- (5);

    \draw[->] (5) -- (9);
    \draw[->] (9) -- (10);
    \draw[->] (10) -- (11);
    \draw[->] (11) -- (1);

\end{tikzpicture}
\caption{The graph $\Theta(a,b,c)$.}
\label{fig:theta}
\end{minipage}
\end{figure}

\begin{definition}
Let $D$ be a directed graph. A {\em subgraph} $D'$ of the graph $D$ is a pair of sets $V(D') \subseteq V(D)$ and $E(D') \subseteq E(D)$ forming a graph. This means that both the initial and the terminal vertices of any arc in $E(D')$ belong to~$V(D')$.
\end{definition}

\begin{lemma}\label{lemma:subgraph1}
Let $D$ be a graph with the largest chainable index of the first type, and suppose that there exists a subgraph $D' \subseteq D$ which is isomorphic to one of the graphs $K(a,b,c)$, $\Theta(a,b,c)$, or~$\infty(a,b)$. Then $D' = D$.
\end{lemma}
\begin{proof}
By definition, all three families of graphs $K, \Theta, \infty$ have exactly one vertex~$w$ with $d^+(w) = 2$, exactly one vertex~$z$ with $d^-(z) = 2$, with all other vertices satisfying $d^+ = d^- = 1$. By Definition~\ref{definition:types}, this means that the vertices of the set $V(D')$ have no incident arcs in the graph $D$ other than the arcs from~$E(D')$. Since Proposition~\ref{proposition:pcisweak} implies that $D$ is connected, we must have $V(D') = V(D)$. It follows from the definitions of the graphs $K$, $\Theta$, and~$\infty$ that a graph $D'$ isomorphic to one of them satisfies $|E(D')| = |V(D')| + 1$. By the definition of a graph of the first type, we obtain $|E(D)| = |V(D)|+1$. Therefore, $E(D') = E(D)$ and $D' = D$.
\end{proof}
\begin{lemma}\label{lemmaallinout}
Let $D$ be a graph with the largest chainable index of the first type, and let $d^+(w) = 2$, $d^-(z) = 2$. Then for any vertex $x\in V(D)$ there exist walks $w \xrightarrow{t_1} x$ and $x \xrightarrow{t_2} z$ for some $t_1,t_2 \in \mathbb{N}$.
\end{lemma}
\begin{proof}
Let $x$ be an arbitrary vertex. The graph $D$ is potentially chainable, hence, by Lemma~\ref{lemma:CIgraph}, for some $t \in \mathbb{N}$ and some vertex $y \neq x$ we have $x \wedge^t y$, that is, there exists a vertex~$v$ such that $x \xrightarrow{t} v$, $y \xrightarrow{t} v$. Let $t_2 \in \mathbb{N}$ be the smallest number such that $x \xrightarrow{t_2} v_0$ and $y \xrightarrow{t_2} v_0$ for some vertex $v_0$. Then $d^-(v_0) \geq 2$. But $z$ is the only vertex in $D$ with $d^- \geq 2$. Hence $v_0 = z$ and $x \xrightarrow{t_2} z$.

The graph $D^T$ is also potentially chainable, and, by Proposition~\ref{propbasic}, in this graph we have $d^-(w) = 2$. Applying the argument of the previous paragraph to the graph $D^T$, we obtain that there is a walk $x\xrightarrow{t_1} w$ in $D^T$. Hence, in the graph $D$ there is a walk $w \xrightarrow{t_1} x$.
\end{proof}

\begin{theorem}\label{THfirsttypeclass}
Let $D$ be a graph with the largest chainable index of the first type. Then $D$ is isomorphic to one of the graphs $K(a,b, c)$, $\infty(a, b)$, or $\Theta(a,b,c)$.
\end{theorem}

\begin{proof}
Let $D$ be a graph with the largest chainable index of the first type. We have a unique vertex $z$ with $d^- = 2$ and a unique vertex $w$ with $d^+ = 2$. By Lemma~\ref{lemma:subgraph1}, it suffices to find a subgraph isomorphic to one of the graphs $K(a,b,c)$, $\infty(a,b)$, or $\Theta(a,b,c)$.

\medskip

First we show that if $w = z$, then $D = \infty(a,b)$. We have $d^+(w) = d^-(w) = 2$. If there is a loop $w \rightarrow w$, then there is exactly one vertex $v \neq w$ such that $w \rightarrow v$. For this vertex, Lemma~\ref{lemmaallinout} holds, that is, $w \rightarrow v \xrightarrow{t} w$ for some $t \in \mathbb{N}$. Let $w \rightarrow v \xrightarrow{b-1} w$ be the shortest path of this form. Then all its internal vertices have $d^+ = d^- = 1$. We obtain a subgraph $\infty(1,b)$, and then apply Lemma~\ref{lemma:subgraph1}.
If there is no loop $w \rightarrow w$, then there are two distinct vertices $v_1, v_2$, different from $w$, such that $w \rightarrow v_1$, $w \rightarrow v_2$. Let the shortest paths from $v_i$ to $w$, $i \in \{ 1, 2\}$, which exist by Lemma~\ref{lemmaallinout}, have the form $v_1 \xrightarrow{a-1} w$, $v_2 \xrightarrow{b-1} w$, $a,b \geq 2$. Then, since $w$ is the only vertex with $d^- > 1$, these paths intersect only at the vertex~$w$. We obtain a subgraph $\infty(a,b)$. By Lemma~\ref{lemma:subgraph1}, $D = \infty(a, b)$.

\medskip

Now let $w\neq z$. By Lemma~\ref{lemmaallinout}, there is a path from $w$ to~$z$. Let $w \xrightarrow{b} z$ be a shortest such path.

If there is no path from $z$ to $w$, then the graph is not strongly connected. By Lemma~\ref{lemmaallinout}, there is a closed walk of positive length from $w$ to $w$, as well as one from $z$ to $z$. Let the shortest such walks, which are simple cycles, have the form $w \xrightarrow{a} w$ and $z \xrightarrow{c} z$, respectively. Then, by the condition on vertex degrees in a graph of the first type, all three walks described above can intersect only at the vertices $z$ or $w$. However, the cycles $w \xrightarrow{a} w$ and $z \xrightarrow{c} z$ cannot intersect, since otherwise there would be a path $z \xrightarrow{t} w$. Moreover, in the path $w \xrightarrow{b} z$ all internal vertices are different from $w$ and~$z$. We obtain a subgraph $K(a,b,c)$, where $w = 1$, $z= a+b$. Then, by Lemma~\ref{lemma:subgraph1}, $D = K(a,b,c)$.

Now suppose that there is a path from $z$ to $w$, and the shortest of them has the form $z \xrightarrow{c} w$. Since $d^+(w) = 2$, the vertex $w$ has one more outgoing arc $w \rightarrow u$, besides the one contained in $w \xrightarrow{b} z$, and $d^-(w) = 1$ implies $u \neq w$. By Lemma~\ref{lemmaallinout}, there is a path from $u$ to $z$, and let $u \xrightarrow{a-1} z$ be a shortest such path. We thus obtain another walk $w \rightarrow u \xrightarrow{a-1} z$, where all internal vertices have $d^- = 1$. We also note that $w$ cannot be an internal vertex of this walk, since $d^-(w) = 1$ and there exists a path $z \xrightarrow{c} w$, where all internal vertices have $d^+ = d^- = 1$. Therefore, $w \rightarrow u \xrightarrow{a-1} z$ is a path, and it cannot intersect the path $w \xrightarrow{b} z$ at internal vertices. We obtain a subgraph~$\Theta(a,b,c)$, where $w = 1$, $z = a+1$. By Lemma~\ref{lemma:subgraph1}, $D = \Theta(a,b,c)$.
\end{proof}

\begin{lemma}\label{lemma:lessPC}
Let $A, B \in \pn$ and $a_{ij} \leq b_{ij}$ for all $i,j$. If the matrix $B$ is not potentially chainable, then $A$ is not potentially chainable either.
\end{lemma}
\begin{proof}
If $a_{ij} \leq b_{ij}$, then $(A^k)_{ij} \leq (B^k)_{ij}$ for all $i, j$ and any~$k \in \mathbb{N}$. Hence, by the definition of a chainable matrix, we obtain the claim of the lemma.
\end{proof}

\begin{definition} \label{IndC}
The {\em cyclicity index} $\sigma(D)$ of a strongly connected directed graph~$D$ is the GCD of the lengths of all its simple directed cycles (if there are no directed cycles, we set $\sigma(D) = 1$). If $D$ is not strongly connected, then $\sigma(D)$ is defined as the LCM of the cyclicity indices of its strongly connected components.
\end{definition}

\begin{theorem}[{\cite[Theorem 3.4.4]{Brualdi}}]\label{theoremsigma}
A strongly connected graph $D$ is primitive if and only if $\sigma(D) = 1$.
\end{theorem}

\begin{definition}
Let $a_1,\ldots, a_k \in \mathbb{N}$ be such that $\gcd(a_1,\ldots, a_k) = 1$. The {\em Frobenius number} $F(a_1,\ldots, a_k)$ is the largest integer $h$ such that the equation
\[
x_1a_1 + \ldots + x_ka_k = h
\]
has no solutions in integers $x_i \geq 0$.
\end{definition}

Frobenius numbers play an important role in various applied problems. For example, the number $F(a_1,\ldots, a_k)$ is the largest monetary amount that cannot be obtained by using only coins of denominations $a_1, \dots, a_k$. Say, using two- and five-cent coins one can obtain any natural number of cents except one and three. Hence $F(2,5) = 3$. The following fact is widely known:
\begin{proposition}\label{prop2frob}
If $\gcd(a,b) = 1$, then $F(a,b) = ab - a - b$.
\end{proposition}

Note that, if $a = 1$ or $b = 1$, then every nonnegative integer $h$ can be expressed as $h = xa + yb$ for some integers $x, y \geq 0$. Hence $F(a,b) = -1$, and the formula in Proposition~\ref{prop2frob} holds true.

\begin{theorem}\label{theorem:typeoneclassification}
A graph $D$ with $n$ vertices and $n+1$ arcs is a graph with the largest chainable index if and only if one of the following three conditions holds:
\begin{enumerate}
    \item $D \cong K(a, b, c)$, where $\gcd(a, c) = 1$.
    \item $D \cong \infty(a,b)$, where $\gcd(a,b) = 1$.
    \item $D \cong \Theta(a,b,c)$, where $\gcd(a + c, b + c) = 1$.
\end{enumerate}
\end{theorem}
\begin{proof}
Theorem~\ref{THfirsttypeclass} shows that a graph with the largest chainable index of the first type must be a graph of the form $K$, $\infty$, or~$\Theta$. The chainable rank of all three graphs equals $n-1$ for any values of the parameters $a,b,c$. By Lemma~\ref{lemmaonlyPC}, it suffices to prove that such graphs are potentially chainable if and only if they satisfy the conditions above.

If a graph is strongly connected, then, by Proposition~\ref{propbasic}(3), its adjacency matrix is irreducible. Therefore, by Theorems~\ref{theoremprimtochain} and~\ref{theoremsigma}, a strongly connected graph is potentially chainable if and only if it is primitive, which is equivalent to the equality $\sigma(D) = 1$. These arguments immediately yield the proof for cases (2) and~(3): the graphs $\infty$ and $\Theta$ are strongly connected, $\sigma(\infty(a,b)) = \gcd(a,b)$, and $\sigma(\Theta(a,b,c)) = \gcd(a + c, b + c)$.

Now consider case (1). We have the graph $K(a,b,c)$, $\gcd(a,c) = d$, which is not strongly connected. Put $k = (a + b + 1) \pmod d$ and construct a graph $K'$ by adding the arc $(a+b, k+1)$ to the graph $K$. This is possible, since $0 \leq k < d \leq a$. In this graph there appears an additional cycle
\[
1 \xrightarrow{b} a+b \xrightarrow{1} 1
\]
of length
\[b + 1 \equiv 0  \pmod d\]
if $k = 0$, and the cycle
\[k+1 \xrightarrow{a-(k+1)} a \xrightarrow{1} 1 \xrightarrow{b} a+b \xrightarrow{1} k+1\]
of length
\[a - k + b + 1 \equiv 0  \pmod d\]
if $k \neq 0$. Together with the cycles
\begin{gather*}
    1\rightarrow 2\rightarrow \ldots \rightarrow a \rightarrow 1,\\
    a + b \rightarrow a + b +1 \rightarrow \ldots \rightarrow a + b + c-1 \rightarrow a + b
\end{gather*}
of lengths $a$ and $c$ it forms the set of all simple directed cycles of the graph $K'$. Hence $\sigma(K') = d$. The resulting graph is strongly connected, so it is potentially chainable if and only if $d = 1$. Since the graph~$K'$ is obtained from $K$ by adding an arc, Lemma~\ref{lemma:lessPC} implies that if $K'$ is not potentially chainable, then $K$ is not potentially chainable either.

We have proved that the graph $K(a,b,c)$ is not potentially chainable for $\gcd(a,c) > 1$. Now we prove that any graph $K(a,b,c)$ with $\gcd(a,c) = 1$ is potentially chainable. Put $l = b + F(a,c) + 1$. Then, by Proposition~\ref{prop2frob}, for any $t \geq l$ we have $1 \xrightarrow{t} a + b$. Moreover, for any vertex~$i$ of the graph $K(a,b,c)$ there is a path $i \xrightarrow{t(i)} a + b$ of length
\[
0 \leq t(i) \leq t_0 = \max\{ a+b-1, c-1\}.
\]
Since $a+ b \xrightarrow{c} a+b$, we obtain that there is a walk $i \xrightarrow{t(i) + mc} a+b$ for any $m \geq 0$. Consequently, for $m_0 = \big\lceil l/c \big\rceil$ we obtain

\[i \xrightarrow{t(i) + m_0c} a+b, \quad 1 \xrightarrow{t(i) + m_0c} a+b.\]

Thus, $1 \wedge^{t(i) + m_0c} i$, whence, by Proposition~\ref{propintersect}, for any vertex~$i$ and any $k \geq t(i) + m_0c$ we have $1 \wedge^{k} i$. Then $1 \wedge^{t_0 + m_0c} i$ for any vertex $i \in V(D)$. By definition, this means that the graph $K(a,b,c)$ is potentially chainable.
\end{proof}

\begin{remark}\label{remark:typeone_n=2}
Theorem~\ref{theorem:typeoneclassification} is valid for graphs with two vertices. However, in this case only $D \cong K(1,1,1)$ and $D \cong \infty(1,2)$ occur among the three families, since the graph $\Theta(a,b,c)$ has at least three vertices. Both graphs are chainable, that is, they have chainable index $1 = n-1$.
\end{remark}

\begin{definition}
\leavevmode
\begin{itemize}
    \item An undirected graph $D$ is called {\em bicyclic} if it contains only two independent cycles.
    \item A directed graph $D$ is called {\em bicyclic} if it contains exactly two directed cycles.
\end{itemize}

\end{definition}

\begin{proposition}[{\cite[p.~583]{WangHuangBelardoMarzi_Theta}}]
If an undirected graph $D$ is connected, then it is bicyclic if and only if $|E(D)| =  |V(D)| + 1$.
\end{proposition}

The graphs $K, \infty$, and $\Theta$ are bicyclic. But not all connected bicyclic directed graphs have this form: some vertices of such graphs may have zero in-degree or out-degree in general. Nevertheless, primitive bicyclic directed graphs always have the form $\Theta(a,b,c)$ or $\infty(a,b)$, and their parameters satisfy the restrictions of Theorem~\ref{theorem:typeoneclassification}. Consequently, such graphs have the largest chainable index. However, the graph $K(a,b,c)$ shows that not all bicyclic graphs with the largest chainable index are strongly connected. The three families of bicyclic graphs $K, \infty, \Theta$ are extremal for a number of invariants and are actively studied; see, e.g., \cite{LinShu, WangHuangBelardoMarzi_Theta, WangHuangBelardoMarzi_Dumbbell}.

The following section shows that not all graphs with the largest chainable index are bicyclic.

\section{Graphs with the largest chainable index of the second type}

Recall that graphs of the second type have exactly two vertices $w_1, w_2$ with $d^+ = 2$, and exactly two vertices $z_1, z_2$ with $d^- = 2$, see Definition~\ref{definition:types}.

\begin{definition}
The graph $W_n'$ with $n \geq 3$ vertices and $n+2$ arcs is obtained from the Wielandt graph $W_n$ by adding the loop $1\rightarrow 1$, see Fig.~\ref{fig:Wielandt-prime}.
\end{definition}

\begin{definition}
The graph $\overline{\infty}(a,b)$, where $a,b \geq 2$, is defined as the graph $K(a,1,b)$ with an added arc $(a + b, 2)$, see Fig.~\ref{fig:two_pentagons_rotated}.
\end{definition}

\begin{figure}[ht]
\centering
\begin{minipage}[b]{0.33\textwidth}
\centering
\begin{tikzpicture}[
    > = Stealth,          
    vertex/.style = {
        font=\sffamily
    }
]
    \def\R{2.0}
    \node (1) at (90:\R)  {$1$};
    \node (2) at (30:\R)  {$2$};
    \node (3) at (-30:\R) {$3$};
    \node (4) at (-90:\R) {$\ldots$};
    \node (5) at (210:\R) {$n-1$};
    \node (6) at (150:\R) {$n$};

    \draw[->] (1) to [out=45,in=135,looseness=6] (1);
    \draw[->] (1) -- (2);
    \draw[->] (2) -- (3);
    \draw[->] (3) -- (4);
    \draw[->] (4) -- (5);
    \draw[->] (5) -- (6);
    \draw[->] (6) -- (1);
    \draw[->] (6) -- (2);
\end{tikzpicture}
 \caption{The graph $W_n'$.}
 \label{fig:Wielandt-prime}
\end{minipage}%
\begin{minipage}[b]{0.67\textwidth}
\centering
\begin{tikzpicture}[
    > = Stealth,          
    vertex/.style = {
        font=\sffamily
    }
]
        \node (1) at (-0.9,  0.9756) {$1$};
        \node (2) at (-0.9, -0.9756) {$2$};
        \node (3) at (-2.7820, -1.8021) {$3$};
        \node (4) at (-4.2180,  0.0000) {$\dots$};
        \node (5) at (-2.7820,  1.8021) {$a$};

        \node (6)  at (1.1820,  0.9756) {$a+1$};
        \node (7)  at (3.3180,  1.8021) {$a+2$};
        \node (8)  at (4.6000,  0.0000) {$\dots$};
        \node (9)  at (3.3180, -1.8021) {$a+b-1\quad$};
        \node (10) at (1.1820, -0.9756) {$a+b$};

        \draw[->] (1) -- (2);
        \draw[->] (2) -- (3);
        \draw[->] (3) -- (4);
        \draw[->] (4) -- (5);
        \draw[->] (5) -- (1);

        \draw[->] (6) -- (7);
        \draw[->] (7) -- (8);
        \draw[->] (8) -- (9);
        \draw[->] (9) -- (10);
        \draw[->] (10) -- (6);

        \draw[->] (1) -- (6);
        \draw[->] (10) -- (2);

    \end{tikzpicture}
    \caption{The graph $\overline{\infty}(a,b)$.}
    \label{fig:two_pentagons_rotated}
\end{minipage}
\end{figure}

\begin{lemma}\label{lemma:subgraph2}
Let $D$ be a graph with the largest chainable index of the second type, and suppose that there exists a subgraph $D'$ isomorphic to $W_m'$ or~$\overline{\infty}$. Then $D' = D$.
\end{lemma}
\begin{proof}
The proof is completely similar to the proof of Lemma~\ref{lemma:subgraph1}.
The graphs $W_m'$ and~$\overline{\infty}$ have two vertices with $d^+ = 2$ and two vertices with $d^- = 2$, while all other vertices satisfy $d^+ = d^- = 1$. Therefore, no vertex of the subgraph $D'$ is joined by an arc to a vertex outside~$D'$. By Proposition~\ref{proposition:pcisweak}, the graph $D$ is connected, hence $V(D') = V(D)$. By the definition of the graphs $W_n'$ and $\overline{\infty}$, $|E(D')| = |V(D')|+2 = |V(D)|+2 = |E(D)|$, whence $D' = D$.
\end{proof}

\begin{lemma}\label{lemmaconnect}
Let $D$ be a graph with the largest chainable index of the second type, and let $d^+(w_i) = d^-(z_j) = 2$, $i, j \in \{ 1, 2\}$. Then for any $i, j \in \{ 1, 2\}$ we have $w_i \rightarrow z_j$.
\end{lemma}
\begin{proof}
By Corollary~\ref{lemmaonewedge}, there is exactly one pair of distinct vertices $u,v$ such that $u \wedge v$. Therefore, for these $u, v$ we have
\[u \rightarrow z_1, \quad u \rightarrow z_2, \quad v\rightarrow z_1, \quad v\rightarrow z_2,\]
whence $u = w_1$, $v = w_2$.
\end{proof}

Consequently, if $u \rightarrow z_i$, then $u = w_1$ or $u = w_2$, and if $w_j \rightarrow v$, then $v = z_1$ or $v = z_2$.

\begin{corollary}\label{corolzwpath}
Let $D$ be a graph with the largest chainable index of the second type, and let $d^+(w_i) = d^-(z_j) = 2$, $i, j \in \{ 1, 2\}$. Then for any vertex $x$ there exist $i, j \in \{ 1, 2\}$ and $t,r \in \mathbb{N}_0$ such that $z_j \xrightarrow{t} x$ and $x \xrightarrow{r} w_i$, and all internal vertices of these paths satisfy $d^+ = d^- = 1$.
\end{corollary}
\begin{proof}
Choose any vertex $x$. Since $D$ is a potentially chainable graph, for some $l \in \mathbb{N}$ and some vertex $y \neq x$ we have $x \wedge^l y$, that is, there is a vertex~$v$ such that $x \xrightarrow{l} v$ and $y\xrightarrow{l} v$. Let $l_0 \in \mathbb{N}$ be the smallest number such that $x \xrightarrow{l_0} v_0$ and $y \xrightarrow{l_0} v_0$ for some vertex $v_0$. Then $d^-(v_0) > 1$. Since $z_1, z_2$ are the only vertices with $d^- > 1$, we obtain $x \xrightarrow{l_0} z_j$. Without loss of generality we may assume that this is the shortest walk of positive length from $x$ to $z_{1}$ or~$z_2$. Lemma~\ref{lemmaconnect} implies that $x \xrightarrow{l_0-1} w_i$, and all internal vertices of this path are different from~$z_{1, 2}$. Moreover, if one of the internal vertices equals $w_1$ or $w_2$, then, by Lemma~\ref{lemmaconnect}, the next vertex in this path is $z_1$ or $z_2$, which contradicts the minimality of the walk $x \xrightarrow{l_0} z_j$. Thus all internal vertices of this path are different from $w_{1,2}$ and $z_{1,2}$, that is, satisfy $d^+ = d^- = 1$.

The graph~$D^T$ is also a graph with the largest chainable index of the second type, so we can apply the arguments of the previous paragraph to it. Since $d^-(z_j) = 2$ in the graph $D$, we have $d^+(z_j) = 2$ in $D^T$, so there is a path $x \xrightarrow{t} z_j$ in~$D^T$. Therefore, we obtain a path $z_j \xrightarrow{t} x$ in the original graph $D$.
\end{proof}

\begin{theorem}\label{theorem:typetwoclassification}
A graph $D$ with $n\geq 3$ vertices and $n+2$ arcs is a graph with the largest chainable index if and only if one of the following two conditions holds:
\begin{enumerate}
    \item $D \cong W_n'$,
    \item $D \cong \overline{\infty}(a,b)$, where $\gcd(a,b) = 1$.
\end{enumerate}
\end{theorem}
\begin{proof}
Let $D$ be a graph with the largest chainable index, and let the vertices $w_1, w_2, z_1, z_2 \in V(D)$ be such that $d^+(w_i) = d^-(z_j) = 2$. By Lemma~\ref{lemmaconnect}, $w_i \rightarrow z_j$, $i,j \in \{ 1,2 \}$.

\medskip

Suppose there is a pair $i,j$ such that $z_i = w_j$. After a suitable renumbering of the vertices $z_{1,2}$ and $w_{1,2}$, we may assume that $z_1 = w_1$, and then
\[
z_1\rightarrow z_1, \quad w_2\rightarrow z_1, \quad z_1 \rightarrow z_2, \quad w_2\rightarrow z_2.
\]
Since $d^+(z_1) = d^-(z_1) = 2$, the vertex $z_1$ has no incident arcs other than the three arcs described above. If, moreover, $w_2 = z_2$, then the vertex $w_2 = z_2$ also has no incident arcs other than those above. Since~$n \geq 3$ and, by Proposition~\ref{proposition:pcisweak}, $D$ is connected, we arrive at a contradiction. Hence the vertices $w_2$ and $z_2$ do not coincide and, apart from the arcs above, must have one more incoming and one more outgoing arc, respectively. Thus, for some vertex $u \neq z_{1,2}$ we have $z_2 \rightarrow u$. By Corollary~\ref{corolzwpath}, there is a path $u\xrightarrow{t} w_i$ for some~$i$, where all internal vertices have $d^+ = d^- = 1$. Since $w_1 = z_1$, $d^-(w_1) = 2$, and there are exactly two distinct incoming arcs $z_1\rightarrow z_1$ and $w_2 \rightarrow z_1$, the case $i =1$ is possible only for $u = w_2$. Let the length of the shortest such path be~$t_0 \in \mathbb{N}_0$. We obtain a subgraph $W_{t_0+3}'$, where $w_1 = z_1 = 1$, $z_2 = 2$, $w_2 = t_0+3$, and $t_0+3 = |V(D')|$. By Lemma~\ref{lemma:subgraph2}, $D = W_{t_0+3}' = W_n'$.

The graph $W_n'$ has a loop $1 \rightarrow 1$. Being strongly connected with $\sigma(W_n') = 1$, it is primitive, and hence potentially chainable. In $W_n'$ there is exactly one unordered pair of distinct intersecting vertices, namely, $1$ and $n$, that is, $\crk(W_n') = n-1$ and, by Lemma~\ref{lemmaonlyPC}, $\ci(W_n') = n-1$.

\medskip

If $z_i\neq w_j$ for all $i,j \in \{ 1, 2 \}$, then $z_1 \rightarrow u_1$, $z_2 \rightarrow u_2$, where the vertices $u_{1,2}, z_{1,2}$ are pairwise distinct. Indeed, if $u_i = z_j$ for some $i,j \in \{1,2\}$, then $z_i\rightarrow z_j$, which is impossible, since $w_1 \rightarrow z_j$, $w_2 \rightarrow z_j$, and $d^-(z_j) = 2$. If $u_1 = u_2$, then $d^-(u_1) = 2$, whence it would follow that $u_1 = z_j$. By Corollary~\ref{corolzwpath}, for $i \in \{ 1, 2\}$
there is a path $u_i \xrightarrow{t_i} w_{j(i)}$ of length $t_i \in \mathbb{N}_0$, where all internal vertices have $d^+ = d^- = 1$. We assume that these paths are the shortest possible, and show that by a suitable renumbering of the vertices $w_{1,2}$ one can achieve $j(i) = i$ for $i \in \{ 1, 2\}$. Suppose that $j(1) = j(2) = 1$. All vertices of the paths $u_i \xrightarrow{t_i} w_1$ have $d^- = 1$. Therefore, the equality $t_1 = t_2$ would contradict the condition $u_1 \neq u_2$, and $t_2 > t_1$ implies that either $u_1 = w_1$, or $u_1$ is an internal vertex of $u_2 \xrightarrow{t_2} w_1$. But $z_1 \rightarrow u_1$ and $d^-(u_1) = 1$, so either $u_2 = z_1$, or the vertex~$z_1$ is an internal vertex of this path, a contradiction.

Consider the paths $z_1 \rightarrow u_1 \xrightarrow{t_1} w_1$ and $z_2 \rightarrow u_2 \xrightarrow{t_2} w_2$. They have lengths $a-1$ and $b-1$, where $a = t_1 + 2$, $b = t_2 + 2$. These paths do not intersect: their endpoints $z_{1,2}, w_{1,2}$ are pairwise distinct, so a common vertex would be internal, but for internal vertices $d^+ = d^- = 1$ holds. We obtain a subgraph $\overline{\infty}(a,b)$, where $w_1 = 1$, $w_2 = a+b$, $z_1 =2$, $z_2 = a+1$, and $a,b \geq 2$. By Lemma~\ref{lemma:subgraph2}, $D = \overline{\infty}(a,b)$.

By Theorem~\ref{theoremprimtochain}, the strongly connected graph $\overline{\infty}(a,b)$ is potentially chainable if and only if it is primitive, which is equivalent to the condition $\sigma(\overline{\infty}(a,b)) = 1$ by Theorem~\ref{theoremsigma}. That is, $\overline{\infty}(a,b)$ is potentially chainable if and only if $\gcd(a,b) = 1$. As in the case of~$W_n'$, the graph $\overline{\infty}(a,b)$ has exactly one pair of intersecting vertices, namely, $1$ and $n=a+b$. Consequently, $\crk(\overline{\infty}(a,b)) = n-1$. It then follows from Lemma~\ref{lemmaonlyPC} that $\ci(\overline{\infty}(a,b)) = n-1$ whenever $\gcd(a,b) = 1$.
\end{proof}

For the reader's convenience and clarity of exposition, we state separately the final complete classification of graphs with the largest chainable index.
\begin{theorem}\label{theorem:fullclassification}
\leavevmode
\begin{enumerate}[label=\Roman*.]
    \item A directed graph $D$ on two vertices has the largest chainable index if and only if one of the following three conditions holds:
    \begin{enumerate}
        \item $D \cong K(1,1,1)$;
        \item $D \cong \infty(1,2)$;
        \item $D$ is the complete directed graph on two vertices, that is,
        \[D = D\left(\begin{pmatrix}
            1 & 1 \\
            1 & 1
        \end{pmatrix}\right).\]
    \end{enumerate}
    \item A directed graph $D$ on $n\geq 3$ vertices has the largest chainable index if and only if one of the following five conditions holds:
    \begin{enumerate}
        \item $D\cong K(a,b,c)$, where $a,b,c \geq 1$ and $\gcd(a,c) = 1$;
        \item $D \cong \infty(a,b)$, where $a,b \geq 1$, $a +b > 2$, and $\gcd(a,b) = 1$;
        \item $D \cong \Theta(a,b,c)$, where $a,b,c \geq 1$, $a + b > 2$, and $\gcd(a+c, b+c) = 1$;
        \item $D \cong W_n'$;
        \item $D \cong \overline{\infty}(a,b)$, where $a,b \geq 2$ and $\gcd(a,b) = 1$.
    \end{enumerate}
\end{enumerate}
\end{theorem}
\begin{proof}
First consider the case of a graph on two vertices. By Lemma~\ref{219}, a graph with the largest chainable index on two vertices has either $n+1 = 3$ or $n+2 = 4$ arcs. By Remark~\ref{remark:typeone_n=2}, a graph $D$ on two vertices with three arcs has the largest chainable index if and only if either $D \cong K(1,1,1)$ or $D \cong \infty(1,2)$. If $D$ has $4$ distinct arcs, then it is the complete directed graph. All entries of its adjacency matrix are positive. Such a graph has chainable index $1 = n-1$, so the first part of the theorem holds.

The second part of the theorem, which classifies graphs on $n\geq 3$ vertices, is obtained by combining Theorems~\ref{theorem:typeoneclassification} and~\ref{theorem:typetwoclassification}.
\end{proof}

\begin{proposition} \label{proposition:nonisomorphic-families}
    If graphs belong to different families among those listed in Theorem~\ref{theorem:fullclassification}, then they are not isomorphic to each other.
\end{proposition}

\begin{proof}
    The case $n = 2$ can be verified directly. Let $n \geq 3$. Graphs from families (1--3) have $n + 1$ arcs, while graphs from families (4--5) have $n + 2$ arcs. The graph $K(a,b,c)$ is not strongly connected, unlike the graphs of all other families. In the graph $\infty(a,b)$ the only two simple cycles have exactly one common vertex, while in the graph $\Theta(a,b,c)$ they have at least one common arc. The graph $W'_n$ contains a loop, whereas the graph $\overline{\infty}(a,b)$ has no loops.
\end{proof}

\begin{definition}
A {\em minimum feedback arc set} of a directed graph is a set of arcs of the smallest size whose removal leaves no directed cycles in the graph.
\end{definition}

The cardinality of a minimum feedback arc set is an analogue of the cyclomatic number of undirected graphs; however, in the general case its computation is a more challenging task.

\begin{corollary}
A graph with the largest chainable index on $n \geq 3$ vertices has a minimum feedback arc set of cardinality at most two.
\end{corollary}
\begin{proof}
In the graph $K(a,b,c)$ it suffices to remove the arcs $(a,1)$ and $(a+b+c-1,a+b)$. Consider the graph $\infty(a,b)$. If $a = 1$ or $b = 1$, then it suffices to remove the loop and any other arc. Otherwise it suffices to remove the arcs $(1,2)$ and $(1,a+1)$. In the graph~$\Theta$ it suffices to remove the arc $(a+1,a+b+1)$ if $c \geq 2$, or the arc $(a+1,1)$ if $c = 1$. In the graph $W_n'$ it suffices to remove the arcs $(1,1)$ and $(2,3)$. In the graph $\overline{\infty}$ it suffices to remove the arcs $(a,1)$ and $(a+1,a+2)$.
\end{proof}

\section{Exponents of graphs with the largest chainable index}

By using the classification of graphs with the largest chainable index obtained above, we can compute their exponents directly.

\begin{lemma}\label{lemma:pathlengths}
\leavevmode
\begin{enumerate}
    \item In $\infty(a,b)$ the set of lengths of walks $1 \xrightarrow{t} 1$ is $\{xa+yb \mid x,y \geq 0\}$.
    \item In $\Theta(a,b,c)$ the set of lengths of walks $1 \xrightarrow{t} 1$ is $\{x(a+c) + y(b+c) \mid x,y \geq 0\}$.
    \item In $\overline{\infty}(a,b)$ the sets of lengths of walks $1 \xrightarrow{t} a+1$, $1 \xrightarrow{t} 2$, $a+b \xrightarrow{t} a+1$, $a+b \xrightarrow{t} 2$ coincide and are equal to $\{1 + xa + yb \mid x,y \geq 0\}$.
\end{enumerate}
\end{lemma}
\begin{proof}
\begin{enumerate}
    \item In the graph $\infty(a,b)$ any walk $1 \xrightarrow{t} 1$ decomposes into a concatenation of walks of the form $1 \xrightarrow{a} 1$ and $1 \xrightarrow{b} 1$.
    \item In the graph $\Theta(a,b,c)$ any walk $1 \xrightarrow{t} 1$ decomposes into a concatenation of walks of the form $1 \xrightarrow{a+c} 1$ and $1 \xrightarrow{b+c} 1$.
    \item Consider the graph~$\overline{\infty}(a,b)$. A walk $1 \xrightarrow{t} a+1$ with $t > 1$ has either the form $1 \xrightarrow{t-1} 1 \rightarrow a+1$ or the form $1 \xrightarrow{t-1} a+b \rightarrow a+1$. Both walks $1 \xrightarrow{t-1} 1$ and $1 \xrightarrow{t-1} a+b$ are concatenations of walks of the form $1 \xrightarrow{a} 1$, $1 \xrightarrow{b} a+b$, $a + b\xrightarrow{b} a+b$, and $a+b \xrightarrow{a} 1$. It can be easily seen that any value of the form $xa + yb$ can be obtained as the length of a certain concatenation of such walks. The lengths of the other three walks $1 \xrightarrow{t} 2$, $a+b \xrightarrow{t} a+1$, and $a+b \xrightarrow{t} 2  $ are considered similarly. \qedhere
\end{enumerate}

\end{proof}

\begin{theorem}
The exponents of graphs with the largest chainable index are as follows:
\begin{enumerate}
    \item $\exp(K(a,b,c)) = 0$;
    \item if $\gcd(a,b) = 1$ and $a < b$, then $\exp(\infty(a,b)) = (a + 1)(b - 1)$;
    \item if $\gcd(a+c, b+c) = 1$ and $a < b$, then
    \[
    \exp(\Theta(a,b,c)) = (a+c+1)(b+c-1) -c;
    \]
    \item $\exp(W_n') = 2n - 3$;
    \item if $\gcd(a,b) = 1$ and $a < b$, then $\exp(\overline{\infty}(a,b)) = (a + 1)(b - 1) + 1.$
\end{enumerate}
\end{theorem}
\begin{proof}

\begin{enumerate}
\item The graph $K$ is not strongly connected, hence it is not primitive.

\item
In the graph $\infty(a,b)$, for any vertices $i,j$ we have $i\xrightarrow{t_1} 1 \xrightarrow{t_2} j$, where $0 \leq t_1, t_2 \leq b -1$. By Proposition~\ref{prop2frob} and Lemma~\ref{lemma:pathlengths}, for all $h \geq F(a,b) + 1 = (a-1)(b-1)$ we have $1\xrightarrow{h} 1$. Consequently,
\[
\exp(\infty(a,b)) \leq 2(b - 1) + (a-1)(b-1) =  (a+1)(b-1).
\]
We now prove that $\exp(\infty(a,b)) \geq (a+1)(b-1)$. Consider the vertices $a+1$ and $a+b-1$. Any walk from $a +1$ to $a+b-1$ of length greater than $b-2$ has the form
\[a +1 \xrightarrow{b-1} 1 \xrightarrow{m} 1 \xrightarrow{b - 1} a+b-1,\]
where $m = xa + yb$, $x,y \geq 0$. The existence of a walk of this form with length $(a+1)(b-1) - 1 = ab + b - a -2$
would mean that
\[a + 1 \xrightarrow{b-1} 1 \xrightarrow{ab - b - a} 1 \xrightarrow{b-1} a+b-1,\]
which contradicts Proposition~\ref{prop2frob} by Lemma~\ref{lemma:pathlengths}.

\item In the graph $\Theta(a,b,c)$, for any vertex $i$ we have $i \xrightarrow{t_1} 1$, where $0 \leq t_1 \leq b + c -1$. By Lemma~\ref{lemma:pathlengths}, for any $x,y \geq 0$ we have
\[1 \xrightarrow{x (a+c) + y(b+c)} 1.\]
Consequently, for any $h \geq F(a+c, b+c) + 1$ there is a walk $i \xrightarrow{t_1} 1 \xrightarrow{h} 1$.
Note that any walk $1 \xrightarrow{h} 1$ ends with the segment
\[
a+1,\, a+b+1,\, \ldots,\, a+b+c-1.
\]
If the vertex $j$ belongs to it, then $1 \xrightarrow{h-l} j$ for some $1 \leq l\leq c$. Consequently, for such $j$ there is a walk $i \xrightarrow{t_1} 1 \xrightarrow{h-l} j$. Then for any $r \geq 0$ there is a walk from $i$ to~$j$ of length $b+c-1 + F(a+c, b+c) + r$.

If $j$ does not belong to the set $\{a+1, a+b+1,\ldots, a+b+c-1\}$, then $1 \xrightarrow{t_2} j$, where $t_2 \leq b-1$. Then for any $h \geq F(a+c, b+c) + 1$ there is a walk $i \xrightarrow{t_1} 1 \xrightarrow{h} 1 \xrightarrow{t_2} j$, that is, for any $r \geq 0$ there is a walk from $i$ to $j$ of length $b + c -1 + F(a+c, b+c) + b + r$. In both cases we obtain that for any $r \geq 0$ there is a walk of length
\begin{align*}
    b + c - 1 + (a+c-1)(b+c-1) - 1 + b + r& \\
    = (a+c+1)(b+c-1) - c + r&,
\end{align*}
that is, $\exp(\Theta(a,b,c)) \leq (a+c+1)(b+c-1) - c$.

Note that, since $a, b \geq 2$ and $a < b$, we have $b \geq 3$. To show that the inequality above becomes an equality, it suffices to consider the walk from $a+2$ to $a+b$:
\[a +2 \xrightarrow{b+c-1} 1 \xrightarrow{h} 1 \xrightarrow{b-1} a+b.\]
If $b+c - 1 + h + b-1 = (a+c+1)(b+c-1) - c - 1$, then
\[
h = (a+c)(b+c) - (a + c) - (b+c).
\]
This contradicts Proposition~\ref{prop2frob}, since, by Lemma~\ref{lemma:pathlengths}, the length of a walk from~$1$ to $1$ can only be equal to $h=x(a+c) + y(b+c)$.

\item In the graph $W'_n$, for any vertices $i,j$ and any $h \geq 0$ there is a walk
\[i \xrightarrow{n-i+1} 1 \xrightarrow{h} 1 \xrightarrow{j-1} j,\]
that is, $i \xrightarrow{n-i+j+h} j$. The expression $n-i+j$ takes a value greater than $2n-3$ only on the pairs $(i,j) \in \{ (1,n), (1,n-1),(2,n)\}$. It can be easily seen that $1 \xrightarrow{j-1+h} j$ for any vertex $j$ and any~$h$. Moreover, there is a walk from $2$ to $n$ of length $2n - 3$, since $2 \xrightarrow{n-2} n$ and $n \xrightarrow{n-1} n$. Thus, for any $i,j$ there is a walk $i \xrightarrow{2n-3} j$, so we obtain $\exp(W_n')\leq 2n-3$.

There is no walk of length $2n-4$ from $2$ to $n$, hence $\exp(W_n') = 2n-3$.

\item In the graph $\overline{\infty}(a,b)$, for any vertex $i$ we have either $i \xrightarrow{t_1} 1$, $0 \leq t_1\leq a-1$, or $i \xrightarrow{t_1} a+b$, $0 \leq t_1 \leq b-1$. For any vertex $j$ we have either $2 \xrightarrow{t_2} j$, $0\leq t_2 \leq a-1$, or $a+1 \xrightarrow{t_2} j$, $0 \leq t_2 \leq b-1$. Using Lemma~\ref{lemma:pathlengths}(3), we obtain that for any vertices $i,j$ and any number $h \geq F(a,b) + 1 = (a-1)(b-1)$ there is a walk
\[i \xrightarrow{t_1+ t_2+h + 1} j, \qquad 0 \leq t_1, t_2 \leq b-1,\]
whence
\[\exp(\overline{\infty}(a,b)) \leq 2(b-1) + (a-1)(b-1) + 1 = (a+1)(b-1)+1.\]

Let us show that the inequality above becomes an equality. Any walk from $a+1$ to $a+b$ has length $xb+b-1$, $x \geq 0$, or $(x+1)a + (y+1)b + b -1$, where $x,y \geq 0$. Thus, the existence of a walk
\[a+1 \xrightarrow{(a+1)(b-1)} a+b\]
is possible in one of two cases:
\begin{enumerate}[label=(\roman*)]
    \item $xb + b-1 = (a+1)(b-1)$. Then $xb = a(b-1)$, which is impossible, since $x \in \mathbb{Z}$ and $b > 1$ is coprime both to $a$ and $b-1$.
    \item $(x+1)a + (y+1)b + b-1 = (a+1)(b-1)$. Then $(x+1)a + yb = ab-a-b$, which contradicts Proposition~\ref{prop2frob}. \qedhere
\end{enumerate}
\end{enumerate}
\end{proof}

\section{The number of graphs on $n$ vertices that attain the largest chainable index}

In this section we obtain formulas for the number of pairwise non-isomorphic graphs with $n\geq 3$ vertices that have the largest chainable index.

The following fact is widely known and is verified directly.

\begin{lemma}
The number of pairs of natural numbers $a,b$ such that $a < b$, $\gcd(a,b) = 1$, and $a + b = n$, where $n \geq 3$, equals $\varphi(n)/2$,
where $\varphi$ is the Euler's totient function.
\end{lemma}

By Proposition~\ref{proposition:nonisomorphic-families}, graphs that belong to different families from Theorem~\ref{theorem:fullclassification} are not isomorphic to each other. We count separately the number of graphs in each of the five families and compute the asymptotic estimates for the obtained values. We do this by counting the number of tuples of parameters $a,b,c$ satisfying the conditions of Theorem~\ref{theorem:fullclassification}.

\begin{enumerate}[widest=Kabcam,leftmargin=*]
    \item[$K(a,b,c)$:] Condition: $2 \leq a+c \leq n$, $\gcd(a,c) = 1$. Note that the additional restriction $a < c$ is not required, since for $a\neq c$ the graph $K(a,b,c)$ is not isomorphic to $K(c,b,a)$. For each $2 \leq m \leq n$ the number of pairs $a, c$ with $a+c = m$ equals $\varphi(m)$, therefore, we have
    \[|K| = \sum_{i=2}^n\varphi(i).\]
    Using the well-known estimate for the sum of the first $n$ values of the Euler's totient function, see~\cite[Theorem~330]{Hardy}, we obtain
    \[|K| = \sum_{i=1}^n \varphi(i) - 1 = \frac{3n^2}{\pi^2} + O(n\log n) - 1 = \frac{3n^2}{\pi^2} + O(n\log n).\]
    \item[$\infty(a,b)$:] Condition: $a + b = n+1$, $\gcd(a,b) = 1$. Here $\infty(a,b)$ is isomorphic to $\infty(b,a)$. Then
    \[|\infty| = \frac{\varphi(n+1)}{2}.\]
    Since $\varphi(n) < n$, we have $|\infty| = O(n)$.
    \item[$\overline{\infty}(a,b)$:] Condition: $a,b > 1$, $a + b = n$, $\gcd(a,b) = 1$. Here $\overline{\infty}(a,b)$ is isomorphic to $\overline{\infty}(b,a)$. Then
    \[|\overline{\infty}|= \frac{\varphi(n)}{2} - 1,\]
    whence $|\overline{\infty}| = O(n)$.
    \item[$W_n'$:] There is one graph up to isomorphism.
    \item[$\Theta(a,b,c)$:]
     Condition: $a + b + c = n+1$, $a< b$, $\gcd(a + c, b + c) = 1$.

     Putting $x = a+c$, $y = b+c$, we obtain the conditions
     \[0 < x < y < n+1, \quad x+y > n+1, \quad \gcd(x,y) = 1. \]
     For the M\"obius function $\mu$ we have the equality
     \[
     \sum_{d \divides \gcd(x,y)} \mu(d) =
     \begin{cases}
         1, & \gcd(x,y)=1,\\
         0, & \text{otherwise}.
     \end{cases}
     \]
     Therefore, the number of pairs $x,y$ satisfying the above conditions equals
     \[|\Theta| = \sum_{ \substack{0 < x < y < n+1, \\  x+y > n+1}} \left( \sum_{d\divides x, d\divides y} \mu(d) \right) = \sum_{d=1}^n \mu(d) N_d,\]
     where $N_d = |M_d|$ is the cardinality of the set
     \[
         M_d = \{(x,y) \mid 0 < x < y < n+1, \; x+y > n+1, \; d\divides x, \; d\divides y\}.
     \]
     For $(x, y) \in M_d$ we have $x = ud$ and $y = vd$. Hence
    \[N_d = \left|\left\{(u,v) \mid 0 < u < v < \tfrac{n+1}{d}, \; u+v > \tfrac{n+1}{d} \right\}\right|.\]
    For a fixed $u$ we obtain
    \[
    \max\left(\tfrac{n+1}{d} -u, u\right) < v < \tfrac{n+1}{d},
    \]
    so the number of pairs $(u,v)$ equals
    \[
    \begin{cases}
        u+\varepsilon, & u < \frac{n+1}{2d},\\
        \frac{n+1}{d} -u + \varepsilon, & u \geq \frac{n+1}{2d},
    \end{cases}
    \]
    where $|\varepsilon| \leq 1$.
    Thus,
    \[
    \left\lvert N_d - \tfrac{n^2}{4d^2} \right\rvert \leq C\tfrac{n}{d},
    \]
    where $C > 0$ is some fixed constant. Since $|\mu(d)| \leq 1$, we have
    \begin{align*}
    \left\lvert |\Theta| - \frac{n^2}{4} \sum_{d=1}^n \frac{\mu(d)}{d^2} \right\rvert &=
    \left\lvert \sum_{d=1}^n \mu(d) N_d - \sum_{d=1}^n\mu(d) \frac{n^2}{4d^2} \right\rvert \\
    &\leq \sum_{d=1}^n |\mu(d)| \left\lvert N_d - \frac{n^2}{4d^2} \right\rvert \\
    &\leq \sum_{d=1}^n C\frac{n}{d} = Cn \sum_{d=1}^n \frac{1}{d}= O(n\log n),
    \end{align*}
    where the last equality is obtained by the Euler--Maclaurin formula.

    Next, we use the equality $\sum_{d\leq n}\mu(d)/d^2 =  6/\pi^2 + O(1/n)$, which follows from the expression for the Dirichlet series $\sum_{d =1}^\infty \mu(d)/d^s =  1/\zeta(s)$ for $\mathop{\mathrm{Re}} s > 1$. Consequently, we obtain
    \[|\Theta| = \frac{n^2}{4}\sum_{d=1}^n \frac{\mu(d)}{d^2} + O(n\log n) = \frac{3n^2}{2\pi^2} + O(n\log n).\]
\end{enumerate}

It follows from the computations that the main contribution to the number of graphs with the largest chainable index is made by graphs of the forms $K$ and~$\Theta$. Summing up, we obtain that the number of graphs with the largest chainable index equals
\begin{align*}
|K| + |\Theta| + |\infty| + |\overline{\infty}| + |W_n'| &= \frac{3n^2}{\pi^2} + \frac{3n^2}{2\pi^2} + O(n\log n) + O(n)\\
&= \frac{9n^2}{2\pi^2} + O(n\log n).
\end{align*}

\bibliographystyle{abbrvurl}
\bibliography{bibfile-arXiv}

\end{document}